\documentclass{amsart}
\usepackage{eucal}
\usepackage{amsthm}

\newcommand{\bbR}{\mathbb{R}}
\newcommand{\eo}{\mathfrak{o}}
\newcommand{\et}{\mathfrak{t}}
\newcommand{\calF}{\mathcal{F}}

\newcommand{\Fsup}{\overline{F}}
\newcommand{\Finf}{\underline{F}}
\newcommand{\Fsupplus}{{\overline{F}}^+\!}
\newcommand{\Finfplus}{{\underline{F}}^+\!}

\newcommand{\abs}[1]{\lvert #1\rvert}
 
\DeclareMathOperator{\vertices}{Vert}
\DeclareMathOperator{\edges}{Edge}
\DeclareMathOperator{\trace}{Tr}

\DeclareMathOperator{\specpoint}{spec_{\text{\rm point}}}
\DeclareMathOperator{\im}{im}

\newtheorem{theorem}{Theorem}[section]
\newtheorem{lemma}[theorem]{Lemma}

\newtheorem{corollary}[theorem]{Corollary}

\theoremstyle{definition}
\newtheorem{definition}[theorem]{Definition}

\theoremstyle{remark}

\numberwithin{equation}{section}

\begin{document}

\title[Approximating spectral invariants of Harper operators II]
{Approximating Spectral invariants of Harper operators on graphs II}
\author{Varghese Mathai}
\address{
Department of Mathematics, University of Adelaide, Adelaide 5005,
Australia}
\email{vmathai@maths.adelaide.edu.au}
\author{Thomas Schick}
\address{FB Mathematik, Universit\"at G\"ottingen, Bunsenstr.~~3, 37073 
G\"ottingen, Germany}
\email{schick@uni-math.gwdg.de} 
\author{Stuart Yates}
\address{Department of Mathematics, University of Adelaide, Adelaide 5005,
Australia}
\email{syates@maths.adelaide.edu.au}

%\date{}
\thanks{%
V.M. and S.Y. acknowledge support from the Australian Research Council.%
}

\begin{abstract}
We study Harper operators and the closely related
discrete magnetic Laplacians (DML) on a graph with a free action
of a discrete group, as defined by Sunada \cite{Sun}.
The spectral density function of the DML is defined using the
von Neumann trace associated with the free action of a
discrete group on a graph. The main result in this paper states
that when the group is amenable, the spectral density function is equal to 
the integrated density of states of the DML that is defined using either
Dirichlet or Neumann boundary conditions. This establishes the main 
conjecture in \cite{MY}. The result is generalized to other self
adjoint operators with finite propagation.
\end{abstract}

\maketitle
\section{Introduction}

In \cite{MY} Mathai and Yates examined the spectral properties of the 
discrete magnetic Laplacian (DML) on a 
graph with a free action of a discrete
group, as defined by Sunada \cite{Sun}. One of the 
results proved in \cite{MY} is that when the group is
amenable and has a finite fundamental domain on the graph, the spectral
density function of the DML --- defined using the von Neumann trace
associated with the discrete group action --- 
can be approximated almost everywhere by the normalized spectral densities
of a series of finite approximations of the DML, consisting of restrictions
of the DML to a regular exhaustion of the graph. The latter approximation
is known in the physics literature as the integrated density of states
of the DML, so that the result can be rephrased as saying that 
the von Neumann spectral density function of the DML is equal to 
the integrated density of states of the DML almost everywhere.
This partial result was known to the experts, cf. \cite{Bel}, \cite{Sh}.

It was shown in
\cite{MY} that in fact the approximation holds at every point
whenever the DML
is associated with a rational weight function, and it was 
conjectured there that the assumption on the weight function was irrelevent.
Based on the idea of the middle named author, the  
conjecture is established here in complete generality, and at the same 
time a much more straightforward argument is presented
for the corresponding result in \cite{MY}. Related arguments can be found
for example in
\cite{DLMSY}, \cite{Eck} and \cite{Elek}, and are known to other
authors, in particular Wolfgang L\"uck and Holger Reich.%%%%

Hereafter, let $X$ be a combinatorial graph with finite fundamental
domain $\calF$ under the free action of an amenable discrete group $\Gamma$.
Each edge $[e]$ of $X$ has associated with it two oriented edges $e$ and
$\overline{e}$ with
opposite orientation. While denoting the set of oriented edges by $\edges X$,
it will be convenient to work with a subset $E^+$ of these edges in which
each combinatorial edge has exactly one oriented representative.
Unless otherwise stated, functions over $X$ will refer to complex-valued
functions over the vertices of $X$.

By F\o{}lner's characterization of amenability (see also \cite{Ad}) there
must exist a \emph{regular exhaustion} of $\Gamma$: a tower of
finite subsets $\Lambda_m\subset\Lambda_{m+1}$, $\bigcup_m \Lambda_m=\Gamma$
satisfying
\begin{equation}
\lim_{k\rightarrow\infty}\;
\frac{\#\partial_\delta \Lambda_k}{\#\Lambda_k}=0\qquad\forall \delta>0,
\end{equation}
where $\partial_\delta \Lambda_k$ denotes the $\delta$-neighborhood
of the boundary of $\Lambda_k$ in $\Gamma$ in the word metric
(with respect to some fixed generating set.)

Let $\{X_m\}$ be a sequence of subgraphs of $X$, with $X_m$ being the
largest subgraph of $X$ contained within the translates by $\Lambda_m$
of the fundamental domain. Then the $X_m$ form a regular exhaustion
of the graph, satisfying
\begin{equation}
\label{e:graphre}
   \lim_{m\rightarrow\infty}\;
   \frac{\#{\vertices \partial_\delta X_m}}{\#{\vertices X_m}}=0
   \quad\forall\delta>0,
\end{equation}
where $\partial_\delta X_m$ refers to the subgraph which is
the $\delta$-neighborhood of the boundary of $X_m$ in $X$ in the
simplicial metric.

The discrete magnetic Laplacian $\Delta_\sigma$ 
is an example of an operator with
bounded propagation, that is there exists a constant $c$ such that 
for any function $f$ with support $Y\subset X$, the function
$\Delta_\sigma f$ has support in the $c$-neighborhood of $Y$.
The following theorem, proved in Section \ref{s:theproof}, provides 
an approximation for the spectral density function of such operators.
This was posed (for the DML) as Conjecture 0.1 of \cite{MY}, where
a proof was obtained for all points $\lambda$ only when the operator
was the DML associated with a rational weight function $\sigma$.

\begin{theorem}[Strong spectral approximation theorem]
  \label{t:ssa}
  Let $X$ be a graph on which there is a free group action by an
  amenable group $\Gamma$, with finite fundamental domain.  Let
  $\big\{X_m\big\}^{\infty}_{m=1}$ be the regular exhaustion of $X$
  corresponding to a regular exhaustion $\Lambda_m$ of $\Gamma$.  Let
  $A$ be a self-adjoint operator of bounded propagation acting on
  $l^2(X)$, which commutes with an adjoint-closed set of twisted
  translation operators, and has spectral density function $F$.
  Construct finite approximations $A_m$ to $A$ by restricting $A$ to
  the space of functions supported on $X_m$, and denote by $F_m$ the
  normalized
  spectral density function of $A_m$,
  \[
  F_m(\lambda)=\frac{1}{\#\Lambda_m}\#\{\text{eigenvalues $\leq \lambda$
    of $A_m$, counting multiplicity}\}.
  \]
  Then
  \begin{equation}
    \lim_{m\to\infty} F_m(\lambda) = F(\lambda),\qquad\forall \lambda\in\bbR.
  \end{equation}
\end{theorem}

\section{The discrete magnetic Laplacian and operators of bounded propagation}
\label{s:obp}

The DML and the Harper operator are defined by a $U(1)$-valued weight
function $\sigma$ on the
edges of the graph $X$ which is \emph{weakly $\Gamma$-invariant}. That is,
for any edge $e$, $\sigma(\overline{e})=\overline{\sigma(e)}$, and for
any $\gamma\in\Gamma$ there exists some $U(1)$-valued function
$s_\gamma$ on the vertices
of $X$ such that
\begin{equation}
\label{e:wge}
\sigma(\gamma e)=\sigma(e)s_\gamma(\et(e))\overline{s_\gamma(\eo(e))}
\end{equation}
where $\eo(e)$ and $\et(e)$ are the origin and terminus of the edge $e$
respectively. The Harper operator $H_\sigma$ is then given by
\[
(H_\sigma f)(v):=
\sum_{\substack{e\in E^+\\ \et(e)=v}} \sigma(e)f(\eo(e))-
\sum_{\substack{e\in E^+\\ \eo(e)=v}} \overline{\sigma(e)}f(\et(e))
\]
and the DML $\Delta_\sigma$ is defined by
\[ 
(\Delta_\sigma f)(v):=
\mathcal{O}(v)f(v)- (H_\sigma f)(v),
\]
where $\mathcal{O}(v)$ is the valence of the vertex $v$.

For some $s_\gamma$ satisfying the weak $\Gamma$-invariance property
for $\sigma$ (\ref{e:wge}), one can define the magnetic translation
operators $T_\gamma$,
\begin{equation}
  \label{e:Tgamma}
  (T_\gamma f)(x)=s_\gamma(\gamma^{-1}x)f(\gamma^{-1}x).
\end{equation}
These form an adjoint-closed set of operators over the space 
$l^2(X)$ of $L^2$ functions on the vertices of $X$, and thus
determine a von Neumann algebra $B(l^2(X))^{T_\Gamma}$ consisting
of all operators over $l^2(X)$ that commute with the $T_\gamma$. This
von Neumann algebra has a finite trace $\trace_{\Gamma}$, defined in
\eqref{e:vNt0}. 
The DML is a self-adjoint operator in this von Neumann algebra, and so one
can define its spectral density function $F$ by
\begin{equation}
F(\lambda)=\trace_{\Gamma} E(\lambda)
\end{equation}
where $E(\lambda)=\chi_{(-\infty,\lambda]}(\Delta_\sigma)$ is the
spectral projection of the DML for the interval $(-\infty,\lambda]$.

More generally one can start with an adjoint-closed set of twisted
translation operators
$T_\gamma$ of the form (\ref{e:Tgamma}),
and then consider a self-adjoint bounded operator $A$
in the von Neumann algebra
$B(l^2(X))^{T_\Gamma}$ which --- like the DML --- has bounded propagation.
\begin{definition}
  An operator $A$ has \emph{propagation bounded by $R$} if for any
  function
  $f\in l^2(X)$, the support of $Af$ is a subset of the $R$-neighborhood
  of the support of $f$ in the simplicial metric.

  The operator $A$ is \emph{weakly $\Gamma$-equivariant} if there
  exists a $\Gamma$-indexed set of functions 
  $t_\gamma\colon  \vertices X\to U(1)$ defining twisted translation operators
  $T_\gamma$: 
  \[
  (T_\gamma f)(x)=t_\gamma(\gamma^{-1}x)f(\gamma^{-1}x),
  \]
  such that $A$ commutes with $T_\gamma$ and $T_\gamma^\ast$ for all
  $\gamma\in\Gamma$. Then $A$ is an element of the associated
  von Neumann algebra $B(l^2(X))^{T_\Gamma}$ with trace
  \begin{equation}
    \label{e:vNt0}
    \trace_{\Gamma}(A) = \sum_{v\in {\mathcal F}} <A\delta_v,\delta_v>.
  \end{equation}
\end{definition}

For each of the subgraphs $X_m$ one can define a finite operator $A_m$
by restricting $A$ to those functions with support on the vertices of
$X_m$; explicitly, $A_m=P_m A i_m$ where $P_m$ is the orthogonal
projection onto $l^2(X_m)$ and $i_m\colon  l^2(X_m)\to l^2(X)$ is the
canonical inclusion.  We define the normalized spectral density
functions of these restricted operators by
\begin{equation}
F_m(\lambda)=\frac{1}{\#\Lambda_m}\#\{\text{eigenvalues $\leq \lambda$
of $A_m$, counting multiplicity}\}.
\end{equation}

Part (i) of Theorem 2.6 of \cite{MY} relates limits of these $F_m$
to the spectral density function $F$ of $A$ when $A$ is the discrete
magnetic Laplacian. The argument deployed though extends trivially to
bounded self-adjoint weakly $\Gamma$-equivariant operators of bounded
propagation generally.

\begin{theorem}[Weak spectral approximation theorem]
  \label{t:wsa}
  Take $X$ to be a graph as above, and 
  let $A$ be a bounded self-adjoint weakly $\Gamma$-equivariant operator
  over $l^2(X)$ of bounded propagation, with restrictions $A_m$ to the
  regular exhaustion $X_m$ of $X$.
  Define the following limits of the normalized spectral density functions,
  \begin{equation}
    \begin{aligned}
      \Fsup(\lambda) &= \limsup_{m\to\infty} F_m(\lambda)
      \\
      \Finf(\lambda) &= \liminf_{m\to\infty} F_m(\lambda)
      \\
      \Fsupplus(\lambda) &= \lim_{\delta\to +0} \Fsup(\lambda+\delta)
      \\
      \Finfplus(\lambda) &= \lim_{\delta\to +0} \Finf(\lambda+\delta).
    \end{aligned}
  \end{equation}
  Then
  \begin{equation}
    F(\lambda)=\Fsupplus(\lambda)=\Finfplus(\lambda),
    \qquad\forall\lambda\in\bbR.
  \end{equation}
\end{theorem}
\begin{proof}
Refer to the proof of Theorem 2.6 of \cite{MY}.
\end{proof}

As a direct consequence, one has that at every point $\lambda$ at which
$F$ is continuous, that is, for every $\lambda\not\in
\specpoint A$,
\begin{equation}
  \label{e:wsa}
  \lim_{m\to\infty} F_m(\lambda)=F(\lambda).
\end{equation}
In particular, equation (\ref{e:wsa}) holds for all but at most
a countable set of points. In the following section a new argument
is presented that extends this result to all $\lambda$.

\section{The strong spectral approximation result}
\label{s:theproof}

A consequence of Theorem \ref{t:wsa} is that if the jump of $F$ at a
discontinuity $\lambda$ can be approximated by the jumps at $\lambda$ of the
piece-wise constant normalized spectral density functions $F_m$, then
the limit of $F_m(\lambda)$ exists and equals $F(\lambda)$.

We first need a small technical lemma.

\begin{lemma}
  \label{l:supinf}
  Given two sequences $\{a_i\}$ and $\{b_i\}$ with $a_i\leq b_i$ for
  all $i$, it follows that $\sup_i a_i \leq \sup_i b_i$ and
  $\inf_i a_i \leq \inf_i b_i$.
  
  In particular, if $\{f_i\}$ is a sequence of monotonically increasing
  functions,
  then $\liminf_{i\to\infty} f_i$ and $\limsup_{i\to\infty} f_i$
  are also monotonically increasing.
\end{lemma}

\begin{proof}
  From $a_i\leq b_i$ for all $i$, it follows that
  \[
  \inf_i a_i \leq a_k \leq b_k \leq \sup_i b_i \qquad\forall k.
  \]
  Then
  $\inf_i a_i \leq b_k$ for all $k$ implies $\inf_i a_i \leq \inf_i b_i$,
  and
  $a_k \leq \sup_i b_i$ for all $k$ implies $\sup_i a_i \leq \sup_i b_i$.

  Consider the sequence $\{f_i\}$ of monotonically increasing functions;
  for any $\delta>0$, $f_i(x+\delta)\geq f_i(x)$ for every $i$. Therefore
  \[
  \inf_{i>k} f_i(x) \leq \inf_{i>k} f_i(x+\delta) \qquad\forall k.
  \]
  Taking the limit or supremum over $k$ then gives
  \[
  \liminf_{k\to\infty} f_k(x) \leq \liminf_{k\to\infty} f_k(x+\delta),
  \]
  that is, the function $\liminf_{k\to\infty} f_k$ is monotonically
  increasing.
  The same argument with $\sup$ instead of $\inf$ gives the corresponding
  inequality for $\limsup _{k\to\infty} f_k$.
\end{proof}

The corollary then follows from Theorem \ref{t:wsa}.

\begin{corollary}
  \label{c:jumps}
  Denote by $D(\lambda)$ the jump in $F$ at $\lambda$,
  \[
  D(\lambda)=\lim_{\delta\to +0} F(\lambda)-F(\lambda-\delta),
  \]
  and similarly let $D_m(\lambda)$ be the jump in the normalized
  spectral density functions $F_m$,
  \[
  D_m(\lambda)=\lim_{\delta\to +0} F_m(\lambda)-F_m(\lambda-\delta).
  \]
  Then for any $\lambda\in\bbR$,
  \begin{equation}
    \lim_{m\to\infty} D_m(\lambda)=D(\lambda)
    \implies
    \lim_{m\to\infty} F_m(\lambda)=F(\lambda).
  \end{equation}
\end{corollary}

\begin{proof}
  Fix some $\epsilon>0$ and $\lambda\in\bbR$, and suppose that
  $\lim_{m\to\infty} D_m(\lambda)=D(\lambda)$. Then we will show that
  $F(\lambda)-3\epsilon<F_m(\lambda)<F(\lambda)+\epsilon$ for all
  sufficiently large $m$.

  By the right continuity and monotonicity of $F$, there exists
  a $\delta>0$ such that $F(x)<F(\lambda)+\frac{\epsilon}{2}$ for
  all $\lambda\leq x<\lambda+\delta$. The $F_m$ approach $F$
  at all points of continuity of $F$; as $F$ has at most a
  countable number of discontinuities there exists some
  $x_0\in[\lambda,\lambda+\delta)$ where
  $\lim_{m\to\infty} F_m(x_0)=F(x_0)<F(\lambda)+\frac{\epsilon}{2}$.
  So for $m$ greater than some $M_1$, it follows that
  $F_m(x_0)<F(\lambda)+\epsilon$. The $F_m$ are monotonically
  increasing, so
  \begin{equation}
    \exists M_1\ \text{such that}\quad
    F_m(\lambda)<F(\lambda)+\epsilon\quad\forall m>M_1.
  \end{equation}

  Now $F(\lambda)=D(\lambda)+\lim_{\delta\to +0} F(\lambda-\delta)$,
  so we can choose $d$ such that
  $F(\lambda-d)+D(\lambda)>F(\lambda)-\epsilon$.
  As the $F_m$ are monotonically increasing functions, $\Finf$ is also
  by Lemma \ref{l:supinf}.
  By Theorem \ref{t:wsa}, $\Finfplus(\mu)=F(\mu)$ for all $\mu$. Therefore
  $\Finf(\mu+\delta)>F(\mu)$ for all $\delta>0$.
  Letting $\mu=\lambda-d$ and $\delta=\frac{d}{2}$ gives
  $\Finf(\lambda-\frac{d}{2})>F(\lambda-d)$ and so for sufficiently large
  $m$, $F_m(\lambda-\frac{d}{2})>F(\lambda-d)-\epsilon$.

  From $F_m(\lambda)\geq F_m(\lambda-\frac{d}{2})+D_m(\lambda)$
  and $F(\lambda-d)> F(\lambda)-D(\lambda)-\epsilon$,
  \begin{equation}
    \begin{aligned}
      \exists M_2\ \text{such that}\quad
      F_m(\lambda) &> F(\lambda-d)-\epsilon+D_m(\lambda)
      \\
      &> F(\lambda) + D_m(\lambda) - D(\lambda) -2\epsilon
      \\
      &\geq F(\lambda)-\abs{D_m(\lambda)-D(\lambda)}-2\epsilon
      \quad\forall m>M_2.
    \end{aligned}
  \end{equation}

  By supposition, $\lim_{m\to\infty} D_m(\lambda)=D(\lambda)$. In
  particular we can pick an $M$ greater than $M_1$ and $M_2$ such that
  $\abs{D_m(\lambda)-D(\lambda)}<\epsilon$ for all $m>M$. Therefore
  \[
  \exists M\ \text{such that}\quad
  m>M \implies F(\lambda)+\epsilon>F_m(\lambda)>F(\lambda)-3\epsilon,
  \]
  and thus
  \[
  \lim_{m\to\infty} D_m(\lambda)=D(\lambda)
  \implies
  \lim_{m\to\infty} F_m(\lambda)=F(\lambda).
  \]
\end{proof}

The following theorem establishes this approximation of the
jumps. Its proof is modelled on the argument in \cite{Elek}.

\begin{theorem}
  \label{t:jumplimit}
  The jump $D(\lambda)$ of $F$ at $\lambda$ is the limit of the jumps
  of the normalized spectral density functions:
  \[
  \lim_{m\to\infty} D_m(\lambda)=D(\lambda)\quad\forall\lambda\in\bbR.
  \]
\end{theorem}
\begin{proof}
  Note that the jumps can be expressed in terms of the dimensions of
  kernels,
  \begin{equation}
    \begin{aligned}
      D(\lambda) &= \dim_\Gamma \ker (A - \lambda), \\
      D_m(\lambda) &= \frac{1}{\#\Lambda_m} \dim \ker (A_m - \lambda).
    \end{aligned}
  \end{equation}

  Consider subsets $Y_m$ of $X_m$ consisting of the $r$-interior of $X_m$,
  where $A$ has propagation bounded by $r$.

  For each $X_m$ define $\dim_{X_k}$ for a subspace $W$ of $l^2(X)$ by
  \[
  \dim_{X_k} W=\frac{1}{\#\Lambda_k}
  \sum_{x\in X_k} \langle P_W \delta_x, \delta_x \rangle
  \]
  where $P_W$ is the orthogonal projection onto $W$. This has the
  following properties for subspaces $W$, $V$ of $l^2(X)$:
  \begin{enumerate}
  \item
    $W \perp V \implies \dim_{X_k} (W \oplus V) = \dim_{X_k} W + \dim_{X_k} V$\!,
  \item
    $W \subset V \implies \dim_{X_k} W \leq \dim_{X_k} V$\!,
  \item
    $P_W \in B(l^2(X))^{T_\Gamma} \implies \dim_{X_k} W = \dim_\Gamma W$\!,
  \item
    $W \subset l^2(X_k)
    \implies \dim_{X_k} W = \frac{1}{\#\Lambda_k} \dim W$\!,
    regarding $l^2(X_k)$ as a subspace of $l^2(X)$.
  \end{enumerate}
  In particular, $\dim_{X_k} l^2(X_k)=\#\calF$ and
  $\dim_{X_k} l^2(Y_k)$ converges to $\#\calF$ as $k$ approaches infinity.

  Let $i_m'\colon  l^2(Y_m)\to l^2(X_m)$ be the inclusion with
  $(i_m' f)(x)=0$ for all $x$ in $X_m\setminus Y_m$.
  Let $A_m'\colon  l^2(Y_m)\to l^2(X_m)$ be the operator $A_m i_m'$. 
  For any function $f$ in $l^2(Y_m)$,
  $(A_m'-\lambda i_m')f=(A_m-\lambda)(i_m'f)$. Therefore as subspaces of
  $l^2(X)$,
  \begin{equation}
    \label{e:ker1}
    \begin{aligned}
      \ker (A_k'-\lambda i_k') &\subset \ker (A_k-\lambda), \\
      \im (A_k'-\lambda i_k') &\subset \im (A_k-\lambda).
    \end{aligned}
  \end{equation}
  Let $D_k'(\lambda)=\frac{1}{\#\Lambda_m}\dim\ker(A_k'-\lambda i_k')$.
  Then taking $\dim_{X_k}$,
  \begin{equation*}
    \begin{aligned}
      D_k'(\lambda) = \dim_{X_k} \ker (A_k'-\lambda i_k')
      &\leq \dim_{X_k} \ker (A_k-\lambda) = D_k(\lambda), \\
      \dim_{X_k} \im (A_k'-\lambda i_k') &\leq \dim_{X_k} \im (A_k-\lambda).
    \end{aligned}
  \end{equation*}
  One also has for any $m$,
  \begin{multline*}
    \lim_{k\to\infty}
    \dim_{X_k} \ker (A_k'-\lambda i_k') +
    \dim_{X_k} \im (A_k'-\lambda i_k')
    =
    \lim_{k\to\infty} \dim_{X_k} l^2(Y_k)
    \\
    = \#\calF
    =
    \dim_{X_m} \ker (A_m-\lambda) +
    \dim_{X_m} \im (A_m-\lambda).
  \end{multline*}
  Consequently,
  \begin{equation}
    \label{e:Dlim1}
    \lim_{k\to\infty}
    D_k(\lambda) -D_k'(\lambda) = 0.
  \end{equation}
  
  Recall that $A_k=P_k A i_k$, where $P_k$ is the projection onto
  $l^2(X_k)$ and $i_k$ is the canonical inclusion of $l^2(X_k)$ into
  $l^2(X)$. In terms of $A$ then, $A_k'=P_k A i_k i_k'$.
  By the bounded propagation of $A$, the support of $A i_k i_k' f$
  is contained within $X_k$ for any $f$ in $l^2(Y_k)$ and hence
  as in (\ref{e:ker1}),
  \begin{equation*}
    \begin{aligned}
      \ker (A_k'-\lambda i_k') &\subset \ker (A-\lambda), \\
      \im (A_k'-\lambda i_k') &\subset \im (A-\lambda).
    \end{aligned}
  \end{equation*}

  The kernel and image of $A-\lambda$ are invariant under the twisted
  translations $T_\gamma$, and so for these subspaces, $\dim_{X_k}$ equals
  $\dim_\Gamma$. It follows that
  \begin{align*}
    D_k'(\lambda) = \dim_{X_k} \ker (A_k'-\lambda i_k')
    &\leq \dim_\Gamma \ker (A_k-\lambda) = D(\lambda),
    \\
    \dim_{X_k} \im (A_k'-\lambda i_k') &\leq \dim_\Gamma \im (A_k-\lambda),
  \end{align*}
  and
  \begin{multline*}
    \lim_{k\to\infty}
    \dim_{X_k} \ker (A_k'-\lambda i_k') +
    \dim_{X_k} \im (A_k'-\lambda i_k')
    \\
    = \#\calF = 
    \dim_\Gamma \ker (A-\lambda) +
    \dim_\Gamma \im (A-\lambda).
  \end{multline*}
  Therefore
  \begin{equation*}
    \lim_{k\to\infty}
    D_k'(\lambda) = D(\lambda),
  \end{equation*}
  which together with Equation (\ref{e:Dlim1}) gives the required limit
  \begin{equation*}
    \lim_{k\to\infty}
    D_k(\lambda) = D(\lambda).
  \end{equation*}
\end{proof}

Theorem \ref{t:ssa} follows from Theorem \ref{t:jumplimit}
and Corollary \ref{c:jumps}.


\begin{thebibliography}{DLMSY}
%{9999}{DLMSY}

\bibitem[Ad]{Ad}
  T.~Adachi,
  A note on the F\o{}lner condition of amenability,
  \textit{Nagoya Math. J.} \textbf{131} (1993) 67--74.

\bibitem[Bel]{Bel} J. Bellissard,
Gap Labeling Theorems for Schr\"{o}dinger's Operators,
 From number theory to
physics (Les Houches, 1989),  538--630, Springer, Berlin, 1992.

\bibitem[DLMSY]{DLMSY} J.~Dodziuk, P.~Linnell, V.~Mathai, T.~Schick
  and S.~Yates, Approximating $L^2$-invariants, and the Atiyah
  conjecture,  \texttt{math.GT/0107049}.

\bibitem[Eck]{Eck} B.~Eckmann, Approximating $\ell_2$-Betti numbers of
  an amenable covering by ordinary Betti numbers, {\em
Comment. Math. Helv.} \textbf{74} (1999) 150-155.

\bibitem[El]{Elek} G.~Elek, On the analytic zero divisor conjecture of
  Linnell, \texttt{math.GR/0111180}.

\bibitem[MY]{MY} V.~Mathai and S.~Yates, Approximating spectral invariants
  of Harper operators on graphs, \textit{J. Functional Analysis,}
   in press, \texttt{math.FA/0006138}.

\bibitem[Sh]{Sh} M. Shubin, Discrete Magnetic Schr\"odinger
operators, {\em
Commun. Math. Phys.},
{\bf 164} (1994), no.2,
259--275.

\bibitem[Sun]{Sun} T.~Sunada,
  A discrete analogue of periodic magnetic Schr\"odinger operators,
  \textit{Contemp. Math.} \textbf{173} (1994), 283--299.

\end{thebibliography}
\end{document}